\documentclass{amsart}
\usepackage{amssymb,latexsym,amsmath,graphicx,graphics,epic,eepic}

\theoremstyle{plain}
\newtheorem{theorem}{Theorem}
\numberwithin{theorem}{section}
\newtheorem{lemma}[theorem]{Lemma}
\newtheorem{proposition}[theorem]{Proposition}
\newtheorem{corollary}[theorem]{Corollary}

\theoremstyle{definition}
\newtheorem{definition}[theorem]{Definition}
\newtheorem{example}[theorem]{Example}

\newtheorem{remark}[theorem]{Remark}


\begin{document}

\title{On Legendrian Surgeries}

\author{Hao Wu}

\address{Department of mathematics and Statistics \\ University of Massachusetts \\ Lederle Graduate Research Tower, 710 North Pleasant Street, Amherst, MA 01003-9305, USA}
\email{wu@math.umass.edu}

\begin{abstract}
We use the Ozsv\'ath-Szab\'o contact invariants to distinguish
between tight contact structures obtained by Legendrian surgeries on
stabilized Legendrian links in tight contact $3$-manifolds. We also
discuss the implication of our result on the tight contact
structures on the Brieskorn homology spheres $-\Sigma(2,3,6n-1)$.
\end{abstract}

\maketitle

\section{Introduction}

A contact structure $\xi$ on an oriented $3$-manifold $M$ is an oriented tangent plane distribution such that there is a
$1$-form $\alpha$ on $M$ satisfying $\xi=\ker\alpha$,
$d\alpha|_{\xi}>0$, and $\alpha\wedge d\alpha>0$. Such a $1$-form is
called a contact form for $\xi$. A curve in $M$ is said to be
Legendrian if it is tangent to $\xi$ everywhere. $\xi$ is said to be
overtwisted if there is an embedded disk $D$ in $M$ such that
$\partial D$ is Legendrian, but $D$ is transversal to $\xi$ along
$\partial D$. A contact structure that is not overtwisted is called
tight.

There are three types of symplectic fillability for contact structures.
\begin{enumerate}
    \item $\xi$ is called Stein fillable if there is a Stein surface $(W,J)$ such that $M=\partial W$ and $\xi=TM\cap J(TM)$.
    \item $\xi$ is called strongly fillable if there is a symplectic $4$-manifold $(W,\omega)$ such that $M=\partial W$, $\omega$ is exact near $M$, and there exists a primitive $\alpha$ of $\omega$ near $M$ satisfying $\xi=\ker(\alpha|_M)$ and $\omega|_\xi>0$.
    \item $\xi$ is called weakly fillable if there is a symplectic $4$-manifold $(W,\omega)$ such that $M=\partial W$ and $\omega|_\xi>0$.
\end{enumerate}

From the works of Eliashberg \cite{E5}, Etnyre and Honda \cite{EH2}, Gromov \cite{Gromov} and Ghiggini \cite{Ghstrongfill,Gh1}, we know
\begin{eqnarray*}
& & \{\text{Stein fillable contact structures}\} \\
& \subsetneq & \{\text{strongly fillable contact structures}\} \\
& \subsetneq & \{\text{weakly fillable contact structures}\} \\
& \subsetneq & \{\text{tight contact structures}\}.
\end{eqnarray*}

The classification problem of overtwisted contact structures was solved by Eliashberg \cite{E1}. The classification of tight contact structures up to isotopy is much more complex, and is only known for limited classes of $3$-manifolds.

Eliashberg \cite{E} and, independently, Weinstein \cite{We} defined the Legendrian surgery, which turns out to be a very useful method of constructing tight contact structures. We will recall Weinstein's construction in details in Section \ref{ssh}. From \cite{E4,EH2,We}, Legendrian surgery is known to preserve the above three types of symplectic fillability. It has been used to produce many interesting examples of tight contact structures.

In many cases, in order to classify tight contact structures, we need to distinguish between tight contact structures constructed by different
Legendrian surgeries. If the Legendrian surgeries are done on the standard contact $S^3$, which is Stain filled by the standard complex $B^4$, then the next two theorems provide an easy criterion.

\begin{theorem}\cite[Theorem 1.2]{LM}\label{stein-dis}
Let $X$ be a smooth $4$-manifold with boundary. Suppose $J_1$,
$J_2$ are two Stein structures with boundary on $X$ with
associated $Spin^\mathbb{c}$-structures $\mathfrak{s}_1$ and
$\mathfrak{s}_2$. If the induced contact structures $\xi_1$ and
$\xi_2$ on $\partial X$ are isotopic, then $\mathfrak{s}_1$ and
$\mathfrak{s}_2$ are isomorphic (and, in particular, have the same
first Chern class).
\end{theorem}

\begin{theorem}\cite[Proposition 2.3]{Go}\label{go-chern}
If $(W,J)$ is obtained from the standard complex $B^4$ by Legendrian surgery on a Legendrian link in the standard contact $S^3$, then the first Chern class $c_1(J)$ of the induced Stein structure $J$ is represented by a cocycle whose value on the $2$-dimensional homology class corresponding to a component of $L$ equals the rotation number of that component.
\end{theorem}

In particular, we have:

\begin{corollary}\label{stein-surgery-dis}
Let $L_1$, $L_2$ be two smoothly isotopic Legendrian links in the standard contact $S^3$ (which is Stein fillable). Suppose that the Thurston-Bennequin numbers of corresponding components of $L_1$ and $L_2$ are equal. Then the Legendrian surgeries on $L_1$, $L_2$ give two tight contact structures $\xi_1$ and $\xi_2$ on the same ambient $3$-manifold. And, if $\xi_1$ and $\xi_2$ are isotopic, then the rotation numbers of corresponding components of $L_1$ and $L_2$ are equal.
\end{corollary}

In practice, we can attain different rotation numbers by stabilizing a Legendrian link in different ways. Then Corollary \ref{stein-surgery-dis} implies that Legendrian surgeries on these stabilized Legendrian links give non-isotopic contact structures. This method can be modified to apply to other Stein fillable contact $3$-manifolds. See, e.g., \cite{GLS,H1,Wu} for applications. The goal of the present paper is to generalize Corollary \ref{stein-surgery-dis} to distinguish between tight contact structures obtained by Legendrian surgeries on stabilized Legendrian links in larger classes of tight contact $3$-manifolds, including all weakly fillable ones. Our main technical tool is the Ozsv\'ath-Szab\'o contact invariant.

\begin{theorem}\label{surgery-main}
Let $(M,\xi)$ be a tight contact $3$-manifold, and
\[
L=K^1\coprod K^2\coprod \cdots \coprod K^m
\]
a Legendrian link in it. For
$j=1,2,\cdots,m$, $i=1,2$, fix integers $s^j$, $p^j_i$, so that
$0\leq p^j_i\leq s^j$. Let $K^j_i$ be the Legendrian knot
constructed from $K^j$ by $p^j_i$ positive stabilizations and
$s^j-p^j_i$ negative stabilizations. Then the Legendrian surgeries
on $L_i=K^1_i\coprod K^2_i\coprod \cdots \coprod K^m_i$ give two
contact structures $\xi_1$ and $\xi_2$ on the same ambient
$3$-manifold $M'$. Assume that $\xi_1$ and $\xi_2$ are isotopic. We
have:

\begin{enumerate}
    \item If $(M,\xi)$ is weakly filled by a symplectic
    $4$-manifold $(W,\omega)$, then, for each $j=1,\cdots,m$,
    \[
      2(p^j_1-p^j_2)\left\{
                   \begin{array}{ll}
                     =0, & \hbox{if $K^j$ represents a torsion element in
    $H_1(W)$;} \\
                     \equiv 0 \mod{d^j}, & \hbox{otherwise, where
    $d^j=\gcd\{\langle\zeta,[K^j]\rangle|\zeta\in H^1(W)\}$.}
                   \end{array}
                 \right.
    \]
    \item If $(M,\xi)$ has non-vanishing Ozsv\'ath-Szab\'o $c^+$-invariant,
then, for each $j=1,\cdots,m$,
    \[
      2(p^j_1-p^j_2)\left\{
                   \begin{array}{ll}
                     =0, & \hbox{if $K^j$ represents a torsion element in
    $H_1(M)$;} \\
                     \equiv 0 \mod{d^j}, & \hbox{otherwise, where
    $d^j=\gcd\{\langle\zeta,[K^j]\rangle|\zeta\in H^1(M)\}$.}
                   \end{array}
                 \right.
    \]
\end{enumerate}
\end{theorem}

The above theorem was proved in the author's attempt to classify
tight contact structures on the Brieskorn homology spheres
$-\Sigma(2,3,6n-1)$. In Section \ref{236}, we will discuss the tight
contact structures on these homology spheres using Theorem
\ref{surgery-main}. It is known to many contact topologists that
there are at most $\frac{n(n-1)}{2}$ tight contact structures on
$-\Sigma(2,3,6n-1)$. Using the tight contact structures on
$M(-\frac{1}{2},\frac{1}{3},\frac{1}{6})$, which are all weakly
fillable, we can give $\frac{n(n-1)}{2}$ different Legendrian
surgery constructions of tight contact structures on
$-\Sigma(2,3,6n-1)$. But it is not known whether these surgeries
give non-isotopic tight contact structures. We will use Theorem
\ref{surgery-main} to show that, among these surgeries, any two
different Legendrian surgeries on the same tight contact structure
on $M(-\frac{1}{2},\frac{1}{3},\frac{1}{6})$ give non-isotopic tight
contact structures on $-\Sigma(2,3,6n-1)$, which implies the following theorem.

\begin{theorem}\label{2n-3}
There are at least $2n-3$ pairwise non-isotopic tight contact
structures on $-\Sigma(2,3,6n-1)$.
\end{theorem}

It is still an open problem whether surgeries on different tight contact structures
on $M(-\frac{1}{2},\frac{1}{3},\frac{1}{6})$ give non-isotopic tight
contact structures on $-\Sigma(2,3,6n-1)$. The author believes that
the answer is yes, and the proof will likely require a better understanding of the Heegaard-Floer homology and the Ozsv\'ath-Szab\'o contact invariants.

\section{Standard symplectic $2$-handle and Legendrian
surgery}\label{ssh}

In this section, we recall Weinstein's construction of the standard
symplectic $2$-handle and the Legendrian surgery in \cite{We}.

Let $(x_1,y_1,x_2,y_2)$ be the standard Cartesian coordinates of
$\mathbb{R}^4$, and
\[
\omega_{st}=dx_1\wedge dy_1 + dx_2\wedge dy_2
\]
the standard symplectic form on $\mathbb{R}^4$. Define
\[
f_2=x_1^2-\frac{y_1^2}{2}+x_2^2-\frac{y_2^2}{2},
\]
\[v_2=\nabla f_2=2x_1\frac{\partial}{\partial
x_1}-y_1\frac{\partial}{\partial y_1}
+2x_2\frac{\partial}{\partial x_2}-y_2\frac{\partial}{\partial
y_2},
\]
and
\[
\alpha_2=\iota_{v_2}\omega_{st}=y_1dx_1+2x_1dy_1+y_2dx_2+2x_2dy_2.
\]
Then $v_2$ is a symplectic vector field, in the sense that
$d(\iota_{v_2}\omega_{st})=\omega_{st}$. Let
\[
X_-=\{(x_1,y_1,x_2,y_2)\in\mathbb{R}^4~|~f_2(x_1,y_1,x_2,y_2)=-1\}.
\]
$X_-$ is positively transverse to $v_2$, and, hence,
$\alpha_2|_{X_-}$ is a contact form. Let
\[
S^1_- = \{(0,y_1,0,y_2)~|~f_2(0,y_1,0,y_2)=-1\}.
\]
This is a Legendrian knot in $(X_-,\alpha_2|_{X_-})$.

\begin{lemma}\cite[Lemma 3.1]{We}\label{handle-construction}
For $A>1$, let
\[
F(x_1,y_1,x_2,y_2) = A(x_1^2+x_2^2)-\frac{y_1^2+y_2^2}{2}-1.
\]
Then the hypersurface
\[
\Sigma = F^{-1}(0)
\]
is positively transverse to $v_2$, and the region
\[
\mathcal{H}_2 = \{(x_1,y_1,x_2,y_2)\in\mathbb{R}^4 ~|~
f_2(x_1,y_1,x_2,y_2)\geq-1, ~F(x_1,y_1,x_2,y_2)\leq0\}
\]
is diffeomorphic to $D^2\times D^2$. Moreover, by choosing
$A\gg1$, we can make $\mathcal{H}_2\cap X_-$ an arbitrarily small
neighborhood of $S^1_-$ in $X_-$.
\end{lemma}

\begin{definition}\label{2-handle}
$(\mathcal{H}_2,\omega_{st}|_{\mathcal{H}_2})$ is called a
standard symplectic $2$-handle.
\end{definition}

\begin{proposition}\cite[Proposition 4.2]{We}\label{isotropic}
Suppose, for $i=1,2$, $(W_i,\omega_i)$ is a symplectic
$4$-manifold, $u_i$ is a symplectic vector field in
$(W_i,\omega_i)$, $M_i$ is a $3$-dimension submanifold of $W$
transverse to $u_i$, and $K_i$ is a Legendrian knot in $M_i$ with
respect to the contact form $\iota_{u_i}\omega_i|_{M_i}$. Then
there is an open neighborhood $U_i$ of $K_i$ in $W_i$, for
$i=1,2$, and a diffeomorphism $\varphi:U_1\rightarrow U_2$, s.t.,
$\varphi^{\ast}(\omega_2|_{U_2})=\omega_1|_{U_1}$,
$\varphi_{\ast}(u_1|_{U_1})=u_2|_{U_2}$, $\varphi(U_1\cap
M_1)=U_2\cap M_2$, $\varphi(K_1)=K_2$.
\end{proposition}

Let $(W,\omega)$ be a symplectic $4$-manifold with boundary, $M$ a
component of $\partial W$, and $\xi$ a contact structure on $M$ so
that $\omega|_{\xi}>0$. Let $K$ be a Legendrian knot in $(M,\xi)$.
By \cite[Lemma 2.4]{EH2}, we isotope $\xi$ near $K$ so that
there exit a neighborhood $U$ of $K$ in $W$, and a non-vanishing
symplectic vector field $v$ defined in $U$, s.t., $v$
transversally points out of $W$ along $U\cap M$, and $\xi|_{U\cap
M} = \ker(\iota_v\omega|_{U\cap M})$. Let $\{\psi_t\}$ be the flow
of $v$. Without loss of generality, we assume there exists
$\tau>0$ such that
\[
U= \bigcup_{0\leq t< \tau}\psi_{-t}(U\cap M).
\]
Choose a small $\varepsilon\in(0,\tau)$. By Proposition
\ref{isotropic}, there is an open neighborhood $V$ of $S^1_-$ in
$\mathbb{R}^4$, and an embedding $\varphi:V\rightarrow U$, s.t.,
$\varphi^{\ast}(\omega)=\omega_{st}$, $\varphi_{\ast}(v_2)=v$,
$\varphi(V\cap X_-)\subset \psi_{-\varepsilon}(U\cap M)$, and
$\varphi(S^1_-)=\psi_{-\varepsilon}(K)$. Choosing $A\gg1$ in Lemma
\ref{handle-construction}, we get a standard symplectic $2$-handle
$\mathcal{H}_2$, such that $\mathcal{H}_2 \cap X_- \subset V$. We
extend the map $\varphi:V \rightarrow U$ by mapping the flow of
$v_2$ to the flow of $v$. Then $\varphi$ becomes a symplectic
diffeomorphism from a neighborhood of $\mathcal{H}_2 \cap X_-$ to
a neighborhood of $K$ in $W$. Now, let
\[
W'=W\cup_\varphi\mathcal{H}_2, ~
\omega' =\left\{%
\begin{array}{ll}
    \omega, & \hbox{on $W$;} \\
    \omega_{st}, & \hbox{on $\mathcal{H}_2$,} \\
\end{array}%
\right., ~\text{and}~
v' =\left\{%
\begin{array}{ll}
    v, & \hbox{on $U$;} \\
    v_2, & \hbox{on $\mathcal{H}_2$.} \\
\end{array}%
\right.
\]
Then $(W',\omega')$ is a symplectic $4$-manifold, and $v'$ is a
symplectic vector field defined in $U\cup_{\varphi}\mathcal{H}_2$,
transversally pointing out of the boundary of $W'$. Let
\[
M'=(M\setminus\mathcal{H}_2)\cup(\mathcal{H}_2\cap\Sigma),
~\text{and}~ \xi'=\left\{%
\begin{array}{ll}
    \xi, & \hbox{on $M\setminus\mathcal{H}_2$;} \\
    \ker\alpha_2, & \hbox{on $\mathcal{H}_2\cap\Sigma$.} \\
\end{array}%
\right.
\]
Then $(M',\xi')$ is the contact $3$-manifold obtained from
$(M,\xi)$ by Legendrian surgery on $K$, and $\omega'|_{\xi'}>0$.

\begin{remark}
If $(M,\xi)$ is weakly fillable, then the above construction gives
$(M',\xi')$ a weak symplectic filling. For a general contact
$3$-manifold $(M,\xi)$, consider the symplectic $4$-manifold
$(M\times I,d(e^t\alpha))$, where $\alpha$ is a contact form for
$\xi$, and $t$ is the variable of $I$. We can carry out the above
construction near $M\times\{1\}$, and get a symplectic cobordism
from $(M,\xi)$ to $(M',\xi')$.
\end{remark}

\section{Ozsv\'ath-Szab\'o invariants and proof of Theorem
\ref{surgery-main}}

Ozsv\'ath and Szab\'o \cite{OS} introduced the
Ozsv\'ath-Szab\'o invariant $c(\xi)$ of a contact structure $\xi$ on
a $3$-manifold $M$. $c(\xi)$ is an element of the quotient
$\widehat{HF}(-M)/\{\pm1\}$ of the Heegaard-Floer homology group of
$-M$, and is invariant under isotopy of $\xi$. $c(\xi)$ vanishes
when $\xi$ is overtwisted. For our purpose, it is more convenient to
use the following variant of the Ozsv\'ath-Szab\'o invariant.

\begin{definition}\cite{Gh1,Pl}\label{c+}
Let $M$ be a closed, oriented $3$-manifold, and
\[
\iota:\widehat{HF}(-M)\rightarrow HF^+(-M)
\]
the canonical map. Define $c^+(\xi)=\iota(c(\xi))$ for any contact structure $\xi$ on
$M$.
\end{definition}

Clearly, $c^+(\xi)$ is also invariant under isotopy of $\xi$, and
vanishes when $\xi$ is overtwisted.

The behavior of Ozsv\'ath-Szab\'o invariants under Legendrian
surgeries is described in the following theorem of Ozsv\'ath and Szab\'o.

\begin{theorem}\cite{OS}\label{OS-surgery}
Let $(M',\xi')$ be the contact $3$-manifold obtained from $(M,\xi)$
by Legendrian surgery on a Legendrian link, then
$F^+_W(c^+(\xi'))=c^+(\xi)$, where $W$ is the cobordism induced by
the surgery.

Specially, this implies that $\xi'$ is tight if $c(\xi)\neq0$.
\end{theorem}

Ghiggini \cite{Gh1} refined Theorem \ref{OS-surgery} to
the following.

\begin{proposition}\cite[Lemma 2.11]{Gh1}\label{Gh1-3.3}
Suppose that $(M',\xi')$ is obtained from $(M,\xi)$ by Legendrian
surgery on a Legendrian link. Then we have
$F^+_{W,\mathfrak{t}}(c^+(\xi'))=c^+(\xi)$, where $W$ is the
cobordism induced by the surgery and $\mathfrak{t}$ is the canonical
$Spin^\mathbb{C}$-structure associated to the symplectic structure
on $W$. Moreover, $F^+_{W,\mathfrak{s}}(c^+(\xi'))=0$ for any
$Spin^\mathbb{C}$-structure $\mathfrak{s}$ on $W$ with
$\mathfrak{s}\neq\mathfrak{t}$.
\end{proposition}

In order to prove Theorem \ref{surgery-main} in the weakly fillable
case, we need to use the Ozsv\'ath-Szab\'o contact invariant twisted by
a $2$-form as defined in \cite{OS1}. Let $(M,\xi)$ be a contact
$3$-manifold with weak symplectic filling $(W,\omega)$, and $B$ an
embedded $4$-ball in the interior of $W$. Consider the element
$\underline{F}^+_{W\setminus B, \mathfrak{s}|_{W\setminus B};
[\omega|_{W\setminus B}]}(c^+(\xi;[\omega|_M]))$ of the group
$\underline{HF}^+(S^3;[\omega|_{S^3}])$, where $S^3=-\partial B$,
$\mathfrak{s}$ is a $Spin^{\mathbb{C}}$-structure on $W$,
$c^+(\xi;[\omega|_M])\in \underline{HF}^+(-M;[\omega|_{M}])$ is the
Ozsv\'ath-Szab\'o contact invariant of $\xi$ twisted by
$[\omega|_M]$, and $\underline{F}^+_{W\setminus B,
\mathfrak{s}|_{W\setminus B}; [\omega|_{W\setminus B}]}$ is the
homomorphism between the two twisted Heegaard-Floer homology groups
induced by the cobordism $W\setminus B$. Note that both
$c^+(\xi;[\omega|_M])$ and $\underline{F}^+_{W\setminus B,
\mathfrak{s}|_{W\setminus B}; [\omega|_{W\setminus B}]}$ are defined
up to an overall multiplication by a factor of the form $\pm T^c$
for some $c\in\mathbb{R}$. To make them absolute, we fix the
auxiliary choices in the constructions of them, including a triple
Heegaard diagram, a base Whitney triangle to define the
homomorphisms, and a representation of $c^+(\xi;[\omega|_M])$. We
also fix a minimal grading generator $\Theta^+$ of $HF^+(S^3)$. Note
that
$\underline{HF}^+(S^3;[\omega|_{S^3}])=HF^+(S^3)\otimes\mathbb{Z}[\mathbb{R}]$.
Let $P_{\xi,\mathfrak{s};[\omega]}\in\mathbb{Z}[\mathbb{R}]$ be the
coefficient of $\Theta^+\otimes1$ in
\[
\underline{F}^+_{W\setminus B,
\mathfrak{s}|_{W\setminus B}; [\omega|_{W\setminus
B}]}(c^+(\xi;[\omega|_M])).
\]
Define a degree on $\mathbb{Z}[\mathbb{R}]$ by setting
$\deg{0}=+\infty$ and $\deg{P}=c_1$ for
\[
P=\sum_{i=1}^m a_i T^{c_i} ~\in~ \mathbb{Z}[\mathbb{R}],
\]
where $a_i\neq0$ and $c_1<\cdots<c_m$. Denote by $\mathfrak{s}_{\omega}$ the canonical $Spin^{\mathbb{C}}$-structure of $(W,\omega)$.

\begin{lemma}\cite[Theorem 4.2]{OS1}\label{min-deg}
\[
\deg P_{\xi,\mathfrak{s}_{\omega};[\omega]}<\deg
P_{\xi,\mathfrak{s};[\omega]}
\]
for any $Spin^{\mathbb{C}}$-structure $\mathfrak{s}$ on $W$ with
$\mathfrak{s}\neq\mathfrak{s}_{\omega}$.
\end{lemma}

\begin{proof}
(Following the proof of \cite[Theorem 4.2]{OS1}.) Fix an open book of $M$ adapted to $\xi$ with connected binding and genus greater than $1$. Eliashberg \cite[Theorem 1.1]{E6} showed than $\omega$ extends over the the Giroux $2$-handle $M\xrightarrow{W_0} M_0$ corresponding to the $0$-surgery on the binding of the open book, where $M_0$ is the surface bundle over $S^1$ resulted from this surgery. Moreover, \cite[Theorem 1.3]{E6} implies that there is a $4$-manifold $V$ with $\partial V = -M_0$, $b_2^+(V)>1$, such that the extension of $\omega$ over $W\cup_M W_0$ further extends to a symplectic structure $\widetilde{\omega}$ on $X=W\cup_M W_0 \cup_{M_0} V$. Let $\widetilde{\mathfrak{s}}_{\widetilde{\omega}}$ be the canonical $Spin^{\mathbb{C}}$-structure of $(X,\widetilde{\omega})$.

Let $\mathfrak{s}$ be any $Spin^{\mathbb{C}}$-structure on $W$ such that $\mathfrak{s}|_M$ is the canonical $Spin^{\mathbb{C}}$-structure of $(M,\xi)$. Using the Composition Law \cite[Theorem 3.9]{OS4} and the arguments in the proof of \cite[Theorem 4.2]{OS1},  one can show that there exists a non-zero element $P\in \mathbb{Z}[\mathbb{R}]$ independent of $\mathfrak{s}$ such that
\begin{eqnarray*}
& & P \cdot \underline{F}^+_{W\setminus B, \mathfrak{s}|_{W\setminus B}; [\omega|_{W\setminus B}]}(c^+(\xi;[\omega|_M])) \\
& = & \sum_{\widetilde{\mathfrak{s}}\in Spin^{\mathbb{C}}(X), ~\widetilde{\mathfrak{s}}|_{W}=\mathfrak{s}, ~\widetilde{\mathfrak{s}}|_{W_0}=\widetilde{\mathfrak{s}}_{\widetilde{\omega}}|_{W_0}, ~\widetilde{\mathfrak{s}}|_{V} = \widetilde{\mathfrak{s}}_{\widetilde{\omega}}|_V} \Phi_{X,\widetilde{\mathfrak{s}}} \cdot T^{\left\langle \omega \cup c_1(\widetilde{\mathfrak{s}}),[X]\right\rangle},
\end{eqnarray*}
where $\Phi_{X,\widetilde{\mathfrak{s}}}$ is the closed $4$-manifold invariant defined in \cite{OS4}. By \cite[Theorem 1.1]{OS5}, the degree of the right hand side of the above equation is equal to $\left\langle \omega\cup c_1(\widetilde{\mathfrak{s}}_{\widetilde{\omega}}),[X]\right\rangle$ if $\mathfrak{s}=\mathfrak{s}_{\omega}$, and is strictly greater than $\left\langle \omega\cup c_1(\widetilde{\mathfrak{s}}_{\widetilde{\omega}}),[X]\right\rangle$ otherwise. This implies the lemma.
\end{proof}

The next two lemmas are technical results needed to prove Theorem
\ref{surgery-main}.

\begin{lemma}\label{bundles}
Let $X$ be a compact manifold with boundary, and $Y$ a closed
submanifold of $X$. Suppose that $L_1$ and $L_2$ are two complex
line bundles over $X$, and there is an isomorphism $\Psi:L_1|_Y
\rightarrow L_2|_Y$. Let $j:(X,\emptyset)\rightarrow(X,Y)$ be the
natural inclusion. Then there exists $\beta~\in~H^2(X,Y)$, such that
$j^{\ast}(\beta)=c_1(L_1)-c_1(L_2)$,  and, for any embedded
$2$-manifold $\Sigma$ in $X$ with $\partial\Sigma\subset Y$, and any
non-vanishing section $v$ of $L_1|_{\partial\Sigma}$, we have
$\langle \beta,[\Sigma] \rangle = \langle c_1(L_1,v),[\Sigma]
\rangle - \langle c_1(L_2,\Psi(v)),[\Sigma] \rangle$, where
$[\Sigma]$ is the relative homology class in $H_2(X,Y)$ represented
by $\Sigma$.
\end{lemma}

\begin{proof}
Denote by $J_i$ the complex structure on $L_i$. Choose a metric
$g_2$ on $L_2|_Y$ compatible with $J_2$, and let
$g_1=\Psi^{\ast}(g_2)$. Consider the complex line bundle
$L=L_1\otimes \overline{L}_2$, where $\overline{L}_2$ is $L_2$ with
the complex structure $-J_2$. Let
$\mathcal{I}:L_2\rightarrow\overline{L}_2$ be the identity map, and
$\overline{\Psi}=\mathcal{I}\circ\Psi$. We define a smooth
non-vanishing section $\eta$ of $L|_Y$ as following: at any point
$p$ on $Y$, pick a unit vector $u_p\in L_1|_p$, and define
$\eta_p=u_p\otimes\overline{\Psi}(u_p)$. It is clear that $\eta_p$
does not depend on the choice of $u_p$ since $\overline{\Psi}$ is
conjugate linear. This gives a smooth non-vanishing section $\eta$
of $L|_Y$. Now, let $\beta=c_1(L,\eta)$. Then
$j^{\ast}(\beta)=c_1(L)=c_1(L_1)-c_1(L_2)$.

Without loss of generality, we assume that $v$ is of unit length.
Choose a section $V_1$ of $L_1|_{\Sigma}$ with only isolated
singularities that extends $v$, and a section $V_2$ of
$L_2|_{\Sigma}$ with only isolated singularities that extends
$\Psi(v)$. Then it is easy to see that
\begin{eqnarray*}
\langle \beta,[\Sigma] \rangle
& = & \text{Sum of indices of
singularities of } ~(V_1 \otimes
\mathcal{I}(V_2)) \\
& = & (\text{Sum of indices of singularities of } ~V_1) \\
& & - (\text{Sum of indices of singularities of } ~V_2) \\
& = & \langle c_1(L_1,v),[\Sigma] \rangle - \langle
c_1(L_2,\Psi(v)),[\Sigma] \rangle.
\end{eqnarray*}
\end{proof}

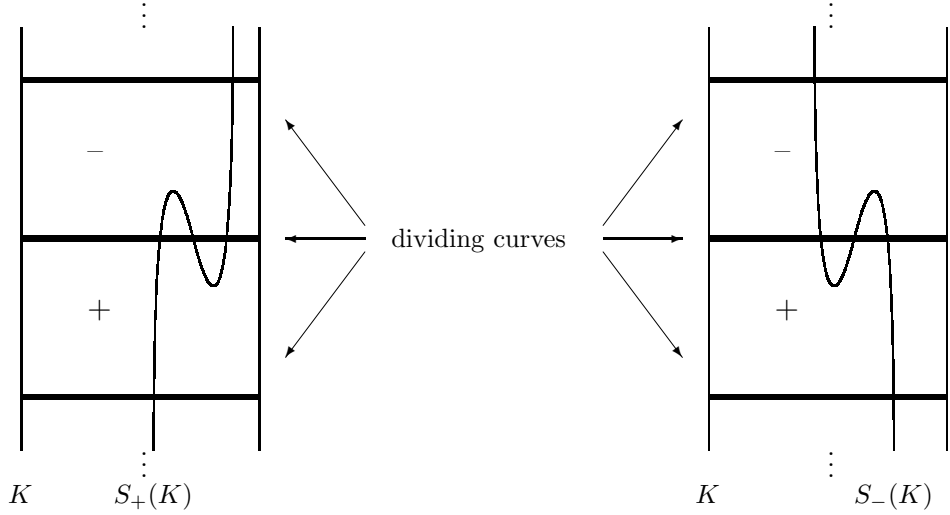
\begin{figure}[ht]

\setlength{\unitlength}{1pt}

\begin{picture}(360,200)(-180,-110)


\linethickness{.5pt}

\put(-175,-80){\line(0,1){160}}

\put(-85,-80){\line(0,1){160}}

\qbezier(-125,-80)(-125,60)(-110,0)

\qbezier(-110,0)(-95,-60)(-95,80)

\put(-130,80){\vdots}

\put(-130,-90){\vdots}

\put(-180,-100){{$K$}}

\put(-140,-100){{$S_+(K)$}}

\put(-150,-30){\Large{+}}

\put(-150,30){\Large{--}}

\linethickness{2pt}

\put(-175,60){\line(1,0){90}}

\put(-175,0){\line(1,0){90}}

\put(-175,-60){\line(1,0){90}}


\linethickness{.5pt}

\put(85,-80){\line(0,1){160}}

\put(175,-80){\line(0,1){160}}

\qbezier(125,80)(125,-60)(140,0)

\qbezier(140,0)(155,60)(155,-80)

\put(130,80){\vdots}

\put(130,-90){\vdots}

\put(80,-100){{$K$}}

\put(140,-100){{$S_-(K)$}}

\put(110,-30){\Large{+}}

\put(110,30){\Large{--}}

\linethickness{2pt}

\put(85,60){\line(1,0){90}}

\put(85,0){\line(1,0){90}}

\put(85,-60){\line(1,0){90}}


\linethickness{.5pt}

\put(-35,-3){dividing curves}

\put(-45,0){\vector(-1,0){30}}

\put(-45,5){\vector(-3,4){30}}

\put(-45,-5){\vector(-3,-4){30}}

\put(45,0){\vector(1,0){30}}

\put(45,5){\vector(3,4){30}}

\put(45,-5){\vector(3,-4){30}}

\end{picture}

  \caption{Positive and Negative
Stabilizations.}\label{stablization-figure}
\end{figure}

Let $K$ be a Legendrian knot in a contact $3$-manifold $(M,\xi)$.
Choose an oriented embedded annulus $\widetilde{A}$ which has $-K$
as one of it is boundary components, and such that the index of the
contact framing of $K$ relative to the framing given by
$\widetilde{A}$ is negative. We can isotope $\widetilde{A}$ relative
to $K$ to make it convex, and such that $K$ has a standard annular
collar $A$ in $\widetilde{A}$. (See, e.g., \cite{H1} for the
definition of standard annular collars.) Then, by Legendrian
Realization Principle \cite[Theorem 3.7]{H1}, we can isotope $A$ relative to $K$ to make
the curved lines in Figure \ref{stablization-figure} Legendrian
without changing the dividing curves. Then these Legendrian curves
are (Legendrianly isotopic to) the positive and negative
stabilizations of $K$. By Giroux's Flexibility, we can again assume
the stabilization has a standard annular collar neighborhood in $A$,
and repeat the above process to obtain repeated stabilizations of
$K$. This observation and \cite[Proposition 4.5]{H1} give:

\begin{lemma}\label{c1-number}
Let $K$ be a Legendrian knot in a contact $3$-manifold $(M,\xi)$.
Then there is an embedded convex annulus $A$ in $M$, such that
$\partial A=(-K)\cup K'$, and $K'$ is (Legendrianly isotopic to) the
repeated stabilization of $K$ obtained by $p$ positive
stabilizations and $s-p$ negative stabilizations. Moreover, if $u$
and $u'$ are the unit tangent vector fields of $K$ and $K'$, then
$\langle c_1(\xi,(-u)\sqcup u'),[A,\partial A]\rangle = 2p-s$.
\end{lemma}

\vspace{.4cm}

\begin{proof}[Proof of Theorem \ref{surgery-main}]  For
notational simplicity, we assume $L=K$ is a Legendrian knot, and
$K_i$, $i=1,2$, is a Legendrian knot obtained from $K$ by $p_i$
positive stabilizations and $s-p_i$ negative stabilizations. The
generalization to Legendrian links is straightforward.

\vspace{.4 cm}

\textbf{Part (1).} We assume that $(M,\xi)$ is weakly filled by
$(W,\omega)$.

First, by \cite[Lemma 2.4]{EH2}, we isotope $\xi$ in near $K$ so
that there is an open neighborhood $U$ of $K$ in $W$ and a
non-vanishing symplectic vector field $v$ defined in $U$, s.t.,
$\xi|_{U\cap M}=\ker(\iota_v\omega|_{U\cap M})$, and $v$
transversally points out of $W$ along $U\cap M$. Let $\{\psi_t\}$ be
the flow of $v$. Without loss of generality, we assume that
$K_i\subset U\cap M$, and
\[
U= \bigcup_{0\leq t< \tau}\psi_{-t}(U\cap M).
\]

Let $W'$ be the smooth $4$-manifold obtained from $W$ by attaching a
$2$-handle to $W$ along $K$ with the framing given by the contact
framing of $K$ plus $s+1$ left twists, and $M'=\partial W'$. Then
the Legendrian surgeries along $K_1$ and $K_2$ give two contact
structures $\xi_1$ and $\xi_2$ on $M'$, and two corresponding
symplectic structures $\omega_1$ and $\omega_2$ on $W'$, such that
$(W',\omega_i)$ is a weak symplectic filling of $(M',\xi_i)$.

\begin{lemma}\label{cohomologous}
We can arrange that $[\omega_1]=[\omega_2]\in H^2(W';\mathbb{R})$.
\end{lemma}
\begin{proof}
Choose a small $\varepsilon\in(0,\tau)$. Let
$N=\psi_{-\varepsilon}(U\cap M)$, and
$\hat{K}_i=\psi_{-\varepsilon}(K_i)$. Then We find a standard
$2$-handle $\mathcal{H}_2$, a neighborhood $V$ of $\mathcal{H}_2\cap
X_-$ in $\mathbb{R}^4$, and an embedding $\varphi_i:V\rightarrow U$,
s.t., $\varphi^{\ast}(\omega)=\omega_{st}$, $\varphi_{\ast}(v_2)=v$,
$\varphi(V\cap X_-)\subset N$, and $\varphi(S^1_-)=\hat{K}_i$. Since
$K_1$ and $K_2$ are isotopic as framed knots,
$\varphi_1|_{\mathcal{H}_2\cap X_-}$ and
$\varphi_2|_{\mathcal{H}_2\cap X_-}$ are isotopic as smooth
embeddings. So there is a smooth isotopy
$\hat{\varphi}_s:\mathcal{H}_2\cap X_-\rightarrow N$, $1\leq
s\leq2$, s.t., $\hat{\varphi}_i=\varphi_i|_{\mathcal{H}_2\cap X_-}$.
After a change of variable in $s$, we assume that
\[
\hat{\varphi}_s=\left\{%
\begin{array}{ll}
    \varphi_1|_{\mathcal{H}_2\cap X_-}, & \hbox{if $1\leq s\leq1.1$;}
\\
    \varphi_2|_{\mathcal{H}_2\cap X_-}, & \hbox{if $1.9\leq s\leq2$.}
\\
\end{array}%
\right.
\]

Let $\widetilde{W}=W\times[1,2]$, and
$\widetilde{\mathcal{H}}_2=\mathcal{H}_2\times[1,2]$. Define
$\widetilde{\omega}$ and $\widetilde{\omega}_{st}$ to be the pull
backs of $\omega$ and $\omega_{st}$ onto $\widetilde{W}$ and
$\widetilde{\mathcal{H}}_2$. And define $\widetilde{v}$ and
$\widetilde{v}_{2}$ to be the lifts of $v$ and $v_2$ to
$U\times[1,2]$ and $\widetilde{\mathcal{H}}_2$ that are tangent to
the horizontal slices $U\times\{s\}$ and $\mathcal{H}_2\times\{s\}$,
$1\leq s\leq2$. Then $\iota_{\widetilde{v}}\widetilde{\omega}$ and
$\iota_{\widetilde{v}_{2}}\widetilde{\omega}_{st}$ are the pull
backs of $\iota_v\omega$ and $\iota_{v_2}\omega_{st}$.

Define $\hat{\Phi}:(\mathcal{H}_2\cap X_-)\times[1,2]\rightarrow
N\times[1,2]$ by $\hat{\Phi}(p,s)=(\hat{\varphi}_s(p),s)$. By
mapping the flow of $\widetilde{v}_{2}$ to the flow of
$\widetilde{v}$, we extend $\hat{\Phi}$ to a diffeomorphism $\Phi$
from a neighborhood of $(\mathcal{H}_2\cap X_-)\times[1,2]$ in
$\widetilde{\mathcal{H}}_2$ to a neighborhood of
$\{\psi_{\varepsilon}\circ\hat{\varphi}_s(S^1_-)~|~1\leq s\leq2\}$
($\subset M\times[1,2]$) in $\widetilde{W}$. Clearly, we have
$\Phi_{\ast}(\widetilde{v}_{2})=\widetilde{v}$, and, near
$(\mathcal{H}_2\cap X_-)\times\{1,2\}$, we have
$\Phi^{\ast}(\widetilde{\omega})=\widetilde{\omega}_{st}$. Consider
the $1$-form $\Phi^{\ast}(\iota_{\widetilde{v}}\widetilde{\omega})$
defined in a neighborhood of $(\mathcal{H}_2\cap X_-)\times[1,2]$ in
$\widetilde{\mathcal{H}}_2$. It equals
$\iota_{\widetilde{v}_{2}}\widetilde{\omega}_{st}$ near
$(\mathcal{H}_2\cap X_-)\times\{1,2\}$. So, there is a $1$-form
$\widetilde{\alpha}$ on $\widetilde{\mathcal{H}}_2$, s.t.,
\[
\widetilde{\alpha}=
\left\{%
\begin{array}{ll}
    \Phi^{\ast}(\iota_{\widetilde{v}}\widetilde{\omega}), & \hbox{near
$(\mathcal{H}_2\cap X_-)\times[1,2]$;} \\
    \iota_{\widetilde{v}_{2}}\widetilde{\omega}_{st}, & \hbox{near
$\mathcal{H}_2\times\{1,2\}$.} \\
\end{array}%
\right.
\]
Define
\[
\widetilde{W}'=\widetilde{W}\cup_{\Phi}\widetilde{\mathcal{H}}_2,
\text{ and }
\widetilde{\omega}'=\left\{%
\begin{array}{ll}
    \widetilde{\omega}, & \hbox{on $\widetilde{W}$;} \\
    d\widetilde{\alpha}, & \hbox{on $\widetilde{\mathcal{H}}_2$.} \\
\end{array}%
\right.
\]
Then $\widetilde{\omega}'$ is a well defined closed $2$-form on
$\widetilde{W}'$.

For $s\in[1,2]$, let
\[
g_s: W' \rightarrow \widetilde{W}'
\]
be the embedding given by $g_s(p)=(p,s)$ for any  point $p$ in $W$
or $\mathcal{H}_2$. Then $\{g_s\}$ is an isotopy of embeddings of
$W'$ into $\widetilde{W}'$. For $i=1,2$, let
$\omega_i=g_i^{\ast}(\widetilde{\omega}')$. Then $(W',\omega_i)$ is
a weak symplectic filling of $(M',\xi_i)$, and
$[\omega_1]=[\omega_2]\in H^2(W';\mathbb{R})$.
\end{proof}

Now we are in the situation that the contact $3$-manifold
$(M',\xi_i)$ is weakly symplectically filled by $(W',\omega_i)$ for
$i=1,2$, $\xi_1$ and $\xi_2$ are isotopic, and
$[\omega_1]=[\omega_2]$. Let $\mathfrak{s}_i$ be the canonical
$Spin^{\mathbb{C}}$-structure of $(W',\omega_i)$. Suppose that
$\mathfrak{s}_1\neq\mathfrak{s}_2$. Then, by Lemma \ref{min-deg}, we
have
\[
\deg P_{\xi_1,\mathfrak{s}_1;[\omega_1]} = \deg
P_{\xi_2,\mathfrak{s}_1;[\omega_2]} > \deg
P_{\xi_2,\mathfrak{s}_2;[\omega_2]},
\]
and, similarly,
\[
\deg P_{\xi_2,\mathfrak{s}_2;[\omega_2]} = \deg
P_{\xi_1,\mathfrak{s}_2;[\omega_1]} > \deg
P_{\xi_1,\mathfrak{s}_1;[\omega_1]}.
\]
This is a contradiction. Thus, $\mathfrak{s}_1=\mathfrak{s}_2$.

Next we construct a symplectic decomposition of $(TW',\omega_i)$ in
a neighborhood of the $2$-handle $\mathcal{H}_2$. First, define a
$2$-plane distribution $\widetilde{\xi}$ on $U$ by
$\widetilde{\xi}|_{\psi_t(p)}=\psi_{t\ast}(\xi_p)$ for $p\in M$. And
let $\widetilde{\eta}=\widetilde{\xi}^{\bot_\omega}$, the
$\omega$-normal bundle of $\widetilde{\xi}$. Clearly, $v$ is a
non-vanishing section of $\widetilde{\eta}$.

Define $\Theta:\mathbb{R}^4\setminus\{0\}\rightarrow Sp(4)$ by
\[
\Theta(x_1,y_1,x_2,y_2) = \frac{1}{\sqrt{4x_1^2+y_1^2+4x_2^2+y_2^2}}
\left(%
\begin{array}{cccc}
   2x_1 & y_1 & -2x_2 & y_2 \\
   -y_1 & 2x_1 & -y_2 & -2x_2 \\
   2x_2 & y_2 & 2x_1 & -y_1 \\
  -y_2 & 2x_2 & y_1 & 2x_1 \\
\end{array}%
\right).
\]
Note that $\Theta$ factors through the natural inclusion of $SU(2)$
into $Sp(4)$ induced by
\[
a+bi \mapsto \left(%
\begin{array}{cc}
  a & -b \\
  b & a \\
\end{array}%
\right).
\]
Since $SU(2)$ is simply connected, we can modify
$\Theta|_{\mathcal{H}_2}$ in a small neighborhood of the
intersection $\mathcal{H}_2\cap\{y_1=y_2=0\}$, and then extend it
into a smooth map $\hat{\Theta}:\mathcal{H}_2\rightarrow Sp(4)$
(c.f. \cite[Proposition 2.3]{Go}). Now let $\{e_1,e_2,e_3,e_4\}$
be the symplectic frame of $T\mathbb{R}^4|_{\mathcal{H}_2}$ defined
by
\[
(e_1,e_2,e_3,e_4) = (\frac{\partial}{\partial x_1},
\frac{\partial}{\partial y_1}, \frac{\partial}{\partial x_2},
\frac{\partial}{\partial y_2})\cdot\hat{\Theta}.
\]

Let $\varphi_i$ be the symplectic attaching map used above to
construct $(W',\omega_i)$, which is a symplectic diffeomorphism from
a neighborhood of $\mathcal{H}_2\cap X_-$ to a neighborhood of $K_i$
in $U$. Note that $\varphi_i$ maps $v_2$
($=\sqrt{4x_1^2+y_1^2+4x_2^2+y_2^2}\cdot e_1$ in the attaching
region) to $v$. So, in the attaching region, $\varphi_i$ identifies
$\widetilde{\xi}$ with the $2$-plane distribution on $\mathcal{H}_2$
spanned by $\{e_3,e_4\}$, and identifies $\widetilde{\eta}$ with the
$2$-plane distribution on $\mathcal{H}_2$ spanned by $\{e_1,e_2\}$.
Let
\[
\widetilde{\xi}_i=\left\{%
\begin{array}{ll}
    \widetilde{\xi}, & \hbox{on $U$;} \\
    \langle e_3,e_4 \rangle, & \hbox{on $\mathcal{H}_2$.} \\
\end{array}%
\right. \text{ and }~
\widetilde{\eta}_i=\left\{%
\begin{array}{ll}
    \widetilde{\eta}, & \hbox{on $U$;} \\
    \langle e_1,e_2 \rangle, & \hbox{on $\mathcal{H}_2$.} \\
\end{array}%
\right.
\]
Then
\[
TW'|_{U\cup_{\varphi_i}\mathcal{H}_2} = \widetilde{\xi}_i \oplus
\widetilde{\eta}_i.
\]
And $\widetilde{\xi}_i$ and $\widetilde{\eta}_i$ are
$\omega_i$-orthogonal to each other. Also, it easy to see that
$\widetilde{\eta}_i$ has a non-vanishing section since we can modify
$v_2$ near the intersection $\mathcal{H}_2\cap\{y_1=y_2=0\}$, and
then extend it to a non-vanishing multiple of $e_1$.

Choose an almost complex structure $J_i$ on
$U\cup_{\varphi_i}\mathcal{H}_2$ compatible with
$\omega_i|_{U\cup_{\varphi_i}\mathcal{H}_2}$ so that
$\widetilde{\xi}_i$ and $\widetilde{\eta}_i$ are complex sub-bundles
of $(TW'|_{U\cup_{\varphi_i}\mathcal{H}_2},J_i)$. Then
$\widetilde{\eta}_i$ becomes a trivial complex line bundle. Note
that $\mathfrak{s}_i|_{U\cup_{\varphi_i}\mathcal{H}_2}$ is the
$Spin^{\mathbb{C}}$-structure associated to $J_i$. There are natural
isomorphisms of complex line bundles
\[
\det (\mathfrak{s}_i)|_{U\cup_{\varphi_i}\mathcal{H}_2} ~\cong~ \det
(TW'|_{U\cup_{\varphi_i}\mathcal{H}_2},J_i) ~\cong~
\widetilde{\xi}_i.
\]
Moreover, there is a natural isomorphism
\[
\det (\mathfrak{s}_i)|_W \cong \det (\mathfrak{s}),
\]
where $\mathfrak{s}$ is the $Spin^{\mathbb{C}}$-structure on $W$
associated to $\omega$.

Let $A_i\subset M$ be the annulus bounded by $(-K)\cup K_i$ given in
Lemma \ref{c1-number}, and
\[
\Sigma_i=A_i\cup(\text{the core of the $2$-handle attached to }K_i),
\]
oriented so that $\partial\Sigma_i=-K$. Then
$[\Sigma_1]=[\Sigma_2]\in H_2(W',W)$. And, by Lemma \ref{bundles},
there exists $\beta~\in~H^2(W',W)$, such that
$j^{\ast}(\beta)=c_1(\det (\mathfrak{s}_1))-c_1(\det
(\mathfrak{s}_2))=0$, and
\begin{eqnarray*}
\langle \beta,[\Sigma_1] \rangle & = & \langle c_1(\det
(\mathfrak{s}_1),-\mu_1),[\Sigma_1] \rangle - \langle c_1(\det
(\mathfrak{s}_2),-\mu_2),[\Sigma_1] \rangle \\
& = & \langle c_1(\widetilde{\xi}_1,-u),[\Sigma_1] \rangle - \langle
c_1(\widetilde{\xi}_2,-u),[\Sigma_1] \rangle \\
& = & \langle c_1(\widetilde{\xi}_1,-u),[\Sigma_1] \rangle - \langle
c_1(\widetilde{\xi}_2,-u),[\Sigma_2] \rangle,
\end{eqnarray*}
where $u$ is the unit tangent vector field of $K$, and $\mu_i$ is
the section of $\det(\mathfrak{s}_i)|_K$ identified with $u$ through
the above isomorphisms.

Denote by $u_i$ the unit tangent vector field of $K_i$. Then $u_i$
extends over the core of the $2$-handle as a non-vanishing multiple
of $e_3$. So, by Lemma \ref{c1-number}, we have $\langle
c_1(\widetilde{\xi}_i,-u),[\Sigma_i]\rangle=\langle
c_1(\xi,(-u)\sqcup u_i),[A_i,\partial A_i]\rangle=2p_i-s$. Thus,
$\langle \beta,[\Sigma_1] \rangle=2(p_1-p_2)$. But, since
$j^{\ast}(\beta)=0$, there exists $\varsigma\in H^1(W)$, s.t.,
$\delta(\varsigma)=\beta$, where $\delta$ is the connecting map in
the long exact sequence of the pair $(W',W)$. So $2(p_1-p_2)=
\langle \delta(\varsigma),[\Sigma_1] \rangle = \langle
\varsigma,-[K] \rangle$. This implies $p_1=p_2$ when $[K]$ is
torsion, and $2p_1\equiv2p_2\mod{d}$ when $[K]$ is non-torsion,
where $d=\gcd\{\langle\zeta,[K]\rangle|\zeta\in H^1(W)\}$.

\vspace{.4 cm}

\textbf{Part (2).} We assume that $c^+(\xi)\neq0$.

Consider the symplectic $4$-manifold $(M\times I,d(e^t\alpha))$,
where $\alpha$ is a contact form for $\xi$, and $t$ is the variable
of $I$. Note that $\frac{\partial}{\partial t}$ is a symplectic
vector field in this setting, and it transversally points out of
$M\times I$ along $M\times \{1\}$. The flow of
$\frac{\partial}{\partial t}$ is the translation in the
$I$-direction. Let $\widetilde{\xi}$ be the $2$-plane distribution
on $M\times I$ generated by translating $\xi$ in the $I$-direction,
and $\widetilde{\eta}=\widetilde{\xi}^{\bot_{d(e^t\alpha)}}$, the
$d(e^t\alpha)$-normal bundle of $\widetilde{\xi}$. Note that
$\frac{\partial}{\partial t}$ is a section of $\widetilde{\eta}$.

We perform Legendrian surgery along $K_i\times \{1\}$. Let
$\varphi_i$ be the symplectic attaching map, which is a symplectic
diffeomorphism from a neighborhood of $S^1_-$ in $\mathcal{H}_2$ to
a neighborhood of $K_i\times \{1\}$ in $M\times I$. Let
\[
W = (M\times I) \cup_{\varphi_1} \mathcal{H}_2 \cong (M\times I)
\cup_{\varphi_2} \mathcal{H}_2.
\]
Then the two Legendrian surgeries give two symplectic structures
$\omega_1$ and $\omega_2$ on $W$, so that $(W,\omega_i)$ is a
symplectic cobordism from $(M,\xi)$ to $(M',\xi_i)$. Similar to the
construction used in Part (1), we construct an $\omega_i$-orthogonal
decomposition
\[
TW = \widetilde{\xi}_i \oplus \widetilde{\eta}_i,
\]
where $\widetilde{\xi}_i|_{M\times I}=\widetilde{\xi}$,
$\widetilde{\eta}_i|_{M\times I}=\widetilde{\eta}$, and, moreover,
$\frac{\partial}{\partial t}$ extends to a non-vanishing section of
$\widetilde{\eta}_i$. Let $J_i$ be an almost complex structure on
$W$ compatible with $\omega_i$ such that both $\widetilde{\xi}_i$
and $\widetilde{\eta}_i$ are complex sub-bundles of $(TW,J_i)$. Then
$\widetilde{\eta}_i$ becomes a trivial complex line bundle over $W$,
and, hence, $c_1(J_i)=c_1(\widetilde{\xi}_i)$.

Let $\mathfrak{s}_i$ be the canonical $Spin^\mathbb{C}$-structure
associated to $J_i$. Then it is also  the canonical
$Spin^\mathbb{C}$-structure associated to $\omega_i$. If
$\mathfrak{s}_1$ and $\mathfrak{s}_2$ are non-isomorphic, according
to Proposition \ref{Gh1-3.3}, we have
\begin{eqnarray*}
F^+_{W,\mathfrak{s}_1}(c^+(\xi_1)) & = &
F^+_{W,\mathfrak{s}_2}(c^+(\xi_2)) ~=~ c^+(\xi) \neq 0 \\
F^+_{W,\mathfrak{s}_1}(c^+(\xi_2)) & = &
F^+_{W,\mathfrak{s}_2}(c^+(\xi_1)) ~=~ 0.
\end{eqnarray*}
But $\xi_1$ and $\xi_2$ are isotopic, this is impossible. So
$\mathfrak{s}_1$ and $\mathfrak{s}_2$ are isomorphic, and, hence,
$c_1(\widetilde{\xi}_1)=c_1(\widetilde{\xi}_2)$.

Let $A_i$  be the annulus in $M\times\{0\}$ bounded by $(-K)\times
\{0\}\cup K_i\times \{0\}$ given by Lemma \ref{c1-number}, and
\[
\Sigma_i=A_i\cup (K_i\times I) \cup(\text{the core of the $2$-handle
attached to }K_i\times \{1\}),
\]
oriented so that $\partial\Sigma=-K\times\{0\}$. Then $\Sigma_1$ and
$\Sigma_2$ are isotopic relative to boundary. And, by Lemma
\ref{bundles}, there exists $\beta~\in~H^2(W,M)$, such that
$j^{\ast}(\beta)=c_1(\widetilde{\xi}_1)-c_1(\widetilde{\xi}_2)=0$,
and
\begin{eqnarray*}
\langle \beta,[\Sigma_1] \rangle & = & \langle
c_1(\widetilde{\xi}_1,-u),[\Sigma_1] \rangle - \langle
c_1(\widetilde{\xi}_2,-u),[\Sigma_1] \rangle \\
& = & \langle c_1(\widetilde{\xi}_1,-u),[\Sigma_1] \rangle - \langle
c_1(\widetilde{\xi}_2,-u),[\Sigma_2] \rangle,
\end{eqnarray*}
where $u$ is the unit tangent vector field of $K\times\{0\}$. Denote
by $u_i$ the unit tangent vector field of $K_i\times \{0\}$. Then,
as in Part (1), $u_i$ extends over $K_i\times I$ and the core of the
$2$-handle without singularities. So, by Lemma \ref{c1-number}, we
have $\langle c_1(\widetilde{\xi}_i,-u),[\Sigma_i]\rangle=\langle
c_1(\xi,(-u)\sqcup u_i),[A_i,\partial A_i]\rangle=2p_i-s$. Thus,
$\langle \beta,[\Sigma_1] \rangle=2(p_1-p_2)$. But, since
$j^{\ast}(\beta)=0$, there exists $\varsigma\in H^1(M)$, s.t.,
$\delta(\varsigma)=\beta$, where $\delta$ is the connecting map in
the long exact sequence of the pair $(W,M\times\{0\})$. So
$2(p_1-p_2)= \langle \delta(\varsigma),[\Sigma_1] \rangle = \langle
\varsigma,-[K] \rangle$. This implies $p_1=p_2$ when $[K]$ is
torsion, and $2p_1\equiv2p_2\mod{d}$ when $[K]$ is non-torsion,
where $d=\gcd\{\langle\zeta,[K]\rangle|\zeta\in H^1(M)\}$.
\end{proof}

\begin{remark}
The weakly fillable case of Theorem \ref{surgery-main} can also be
proved using the monopole invariant defined by Kronheimer and
Mrowka \cite{KM}. Indeed, in Part (1) of the proof, after
proving Lemma \ref{cohomologous}, we are in the situation where
$\xi_1$ and $\xi_2$ are isotopic, and $[\omega_1]=[\omega_2]\in
H^2(W';\mathbb{R})$. After a possible isotopy supported near $M'$,
we assume that $\xi_1=\xi_2=\xi'$. Let $\mathfrak{s}_i\in
Spin^\mathbb{C}(W',\xi')$ be the element associated to $\omega_i$.
Then, by \cite[Theorems 1.1 and 1.2]{KM}, we have
$[\omega_1]\cup(\mathfrak{s}_2-\mathfrak{s}_1)\geq0$, and
$[\omega_2]\cup(\mathfrak{s}_1-\mathfrak{s}_2)\geq0.$ But
$[\omega_1]=[\omega_2]$. Thus,
$[\omega_1]\cup(\mathfrak{s}_2-\mathfrak{s}_1)=0$. And, according
to \cite[Theorem 1.2]{KM}, this implies that
$\mathfrak{s}_1=\mathfrak{s}_2$ as elements of
$Spin^\mathbb{C}(W',\xi')$, and, specially, that
$c_1(\mathfrak{s}_1) = c_1(\mathfrak{s}_2)$. Then we can repeat
the rest of Part (1) of the proof, and prove the weakly fillable
case of the theorem.
\end{remark}

\section{Tight contact structures on Brieskorn homology spheres
$-\Sigma(2,3,6n-1)$}\label{236}

A small Seifert fibered manifold is a  $3$-manifold Seifert fibered
over $S^2$ with $3$ singular fibers. We denote by $M(r_1,r_2,r_3)$
the small Seifert fibered manifold with $3$ singular fibers with
coefficients $r_1$, $r_2$ and $r_3$.

The classification of tight contact structures on a small Seifert
fibered manifold is a hard problem. When the
Euler number of the small Seifert fibered manifold is not $-1$ or
$-2$, these tight contact structures are all Stein
fillable, and are classified in \cite{GLS,Wu}. Note all these manifolds are
$L$-spaces, i.e. have Heegaard-Floer homology like that of a lens
space. There are also partial results when the Euler number is
$-1$ or $-2$, and the manifold is an $L$-space (see e.g.
\cite{GLS2}). In solving these examples, the use of untwisted
Ozsv\'ath-Szab\'o contact invariant is essential. It appears that
the classification is much harder to achieve when the the small
Seifert fibered manifold is not an $L$-space.

Brieskorn homology sphere $-\Sigma(2,3,6n-1)$ is the small Seifert
fibered manifold $M(-\frac{1}{2},\frac{1}{3},\frac{n}{6n-1})$, which is not an $L$-space when
$n\geq2$. These appear to be good examples
of non-$L$-space small Seifert fibered manifolds to start with. In \cite{GS},
Ghiggini and Sch\"{o}nenberger showed that there is a
unique tight contact structure on $-\Sigma(2,3,11)$. This method was
extended to classify contact structures on $-\Sigma(2,3,17)$ in
\cite{Gh1}. Next we discuss the generalization of their method.

Let $\Sigma$ be an oriented three-hole sphere with boundary components $C_1$, $C_2$ and $C_3$.
Then $-\partial\Sigma\times S^1=T_1+T_2+T_3$, where the "$-$" sign
means reversing the orientation and $T_i=-C_i\times S^1$. We identify $T_i$ to
$\mathbb{R}^2/\mathbb{Z}^2$ by identifying $-C_i\times\{\text{pt}\}$ to $(1,0)^T$,
and $\{\text{pt}\}\times S^1$ to $(0,1)^T$. Also, for $i=1,2,3$,
let $V_i=D^2\times S^1$, and identify $\partial V_i$ with
$\mathbb{R}^2/\mathbb{Z}^2$ by identifying a meridian $\partial
D^2 \times \{\text{pt}\}$ with $(1,0)^T$ and a longitude
$\{\text{pt}\}\times S^1$ with $(0,1)^T$.

Define diffeomorphism $\varphi_i: \partial V_i \rightarrow T_i$ by the following matrices.
\[
\varphi_1 =
\left(%
\begin{array}{cc}
  2 & -1 \\
  1 & 0 \\
\end{array}%
\right),
\hspace{.5cm}
\varphi_2 =
\left(%
\begin{array}{cc}
  3 & 1 \\
  -1 & 0 \\
\end{array}%
\right),
\hspace{.5cm}
\varphi_3 =
\left(%
\begin{array}{cc}
  6n-1 & 6 \\
  -n & -1 \\
\end{array}%
\right).
\]
Then
\[
-\Sigma(2,3,6n-1) \cong M(-\frac{1}{2},\frac{1}{3},\frac{n}{6n-1}) \cong (\Sigma
\times S^1)\cup_{(\varphi_1\cup\varphi_2\cup\varphi_3)}(V_1\cup V_2\cup V_3).
\]
Note that each $S^1$-fiber in the product $\Sigma\times S^1$ becomes a regular fiber of the Seifert fibration, and the framing of the $S^1$-fiber from the product is the same as the standard framing of a regular fiber of the Seifert fibration. Also, the core curve of each $V_i$ becomes a singular fiber of the Seifert fibration, and our choice of the longitude of $\partial V_i$ gives each singular fiber a framing. If $\xi$ is a contact structure on $-\Sigma(2,3,6n-1)$, and $K$ is a Legendrian regular fiber (resp. Legendrian singular fiber) of the Seifert fibration, then the twisting number $t(K)$ of $K$ is defined to be the index of the contact framing of $K$ with respect to the standard framing (resp. the framing we chose). We define
\[
t(\xi)=\max\{t(K)~| ~K \text{ is a Legendrian regular fiber.}\}
\]

Etnyre and Honda \cite{EH} showed that $-\Sigma(2,3,5)$ does not admit tight contact structures. So we assume that $n\geq2$ in the discussions below. Next two lemmas are proved following the arguments in \cite[Subsection 4.2]{GS}. Similar methods were also used in e.g. \cite{Wu}.

\begin{lemma}\label{neg-twist}
If $\xi$ is a tight contact structure on $-\Sigma(2,3,6n-1)$, then $t(\xi)\leq -2$.
\end{lemma}

\begin{proof}
We prove the lemma in two steps: first prove that $t(\xi)<0$, and then prove that $t(\xi)\neq-1$.

Assume $t(\xi)\geq0$. Then we can find a Legendrian regular fiber $F$ with twisting number $0$. After possibly an isotopy, assume $F$ is contained in the piece $\Sigma\times S^1$. Let $F_i$ be a Legendrian knot $C^0$-close to the core curve of $V_i$. After repeated stabilization of $F_i$, we assume that $t(F_i)=n_i<<0$. After isotopy, assume that $V_i$ is a standard neighborhood of $F_i$. Then $\partial V_i$ is convex and has two parallel dividing curves of slope $\frac{1}{n_i}$. Now use the coordinates of $T_i$. Then the slopes of dividing curves of $T_1$, $T_2$ and $T_3$ are $s_1=\frac{n_1}{2n_1-1}$, $s_2=-\frac{n_2}{3n_2+1}$ and $s_3=-\frac{nn_3+1}{(6n-1)n_3+6}$, respectively. Since $n_i<<0$, we have $s_1>0$, $s_2>-\frac{1}{2}$ and $s_3>-\frac{1}{5}$. Then we can isotope $T_i$ as in \cite[Subsection 4.2.2]{GS} and get a decomposition
\[
\Sigma \times S^1 = (\Sigma' \times S^1) \cup (T_1 \times [0,1]) \cup (T_2 \times [0,1]) \cup (T_3 \times [0,1]),
\]
such that
\begin{itemize}
    \item $\Sigma'$ is a three-sphere in $\Sigma$ with \[ \partial \Sigma' \times S^1 = (-T_1\times\{1\}) \cup (-T_2\times\{1\}) \cup (-T_3\times\{1\});\]
    \item $\xi|_{T_i \times [0,1]}$ is a minimal twisting tight contact structure with minimal convex boundary;
    \item The slopes of dividing curves on $T_1\times \{0\},T_2\times \{0\},T_3\times \{0\}$ are $0,-\frac{1}{2},-\frac{1}{5}$, respectively, and the slopes of dividing curves on $T_1\times \{1\},T_2\times \{1\},T_3\times \{1\}$ are $\infty$.
\end{itemize}

Then we can follow the arguments in the proof of \cite[Theorem 4.14]{GS} to show that $\xi$ must be overtwisted. This contradiction shows that $t(\xi)<0$.

Now assume that $t(\xi)=-1$. Let $F\subset \Sigma\times S^1$ be a Legendrian regular fiber with $t(F)=-1$, and $V_i$ a standard neighborhood of a Legendrian singular fiber $F_i$ with $t(F_i)=n_i<<0$. For $i=1,2$, connect $F$ to $\partial V_i$ by a vertical convex annulus $A_i$ that intersects the dividing curves of $\partial V_i$ efficiently. By Imbalance Principle \cite[Proposition 3.17]{H1}, there is a $\partial$-parallel dividing curve on $A_i$ along $A_i \cap (\partial V_i)$. Using the bypass from this $\partial$-parallel dividing curve, by the Twisting Number Lemma \cite[Lemma 4.4]{H1}, we can increase $n_i$ by $1$. Repeat this procedure, we can increase $n_1,n_2$ up to $n_1=0$, $n_2=-1$. When measured in the coordinates of $T_i$, the dividing curves on $T_1$ and $T_2$ have slopes $0$ and $-\frac{1}{2}$. Connecting $T_1$ to $T_2$ by a vertical convex annulus $A$ in $\Sigma\times S^1$ with $\partial A$ intersecting the dividing curves of $T_1,T_2$ efficiently. Then, by Imbalance Principle, there is a $\partial$-parallel dividing curve on $A$ along $A\cap T_2$. Adding the bypass from this dividing curve to $T_2$, we change the slope of dividing curves of $T_2$ to $-1$. Connect $T_1$ to this new $T_2$ by a vertical convex annulus $A'$ in $\Sigma\times S^1$ with $\partial A'$ intersecting the dividing curves of $T_1,T_2$ efficiently. If there are $\partial$-parallel dividing curves on $A'$, then, by Legendrian Realization Principle \cite[Theorem 3.7]{H1}, we can find a Legendrian regular fiber with twisting number $0$. This is a contradiction. If there are no $\partial$-parallel dividing curves on $A'$, cut $\Sigma\times S^1$ along $A'$ and smooth the edges. This gives us torus $T$ isotopic to $T_3$ with dividing curves of slope $0$. Note that the slope of dividing curves of $T_3$ is negative since $n_3<<0$. By \cite[Proposition 4.16]{H1}, there is a torus isotopic to $T_3$ with vertical dividing curves (isotopic to a regular fiber.) By the Legendrian Realization Principle, we can again find a Legendrian regular fiber with twisting number $0$, which is a contradiction. This implies that $t(\xi)\neq-1$.
\end{proof}

\begin{lemma}\label{upperbound}
There are at most $\frac{n(n-1)}{2}$ pairwise non-isotopic tight
contact structures on the Brieskorn homology sphere
$-\Sigma(2,3,6n-1)$.
\end{lemma}
\begin{proof}
Let $\xi$ be a tight contact structure on $-\Sigma(2,3,6n-1)$ with $t(\xi)=t$, where $t\leq -2$ by Lemma \ref{neg-twist}. Let $F\subset \Sigma\times S^1$ be a Legendrian regular fiber with $t(F)=t$. Isotope $V_i$ into a standard neighborhood of a Legendrian singular fiber $F_i$ with $t(F_i)=n_i<<0$. For $i=1,2$, connect $F$ to $\partial V_i$ by a vertical convex annulus $A_i$ that intersects the dividing curves of $\partial V_i$ efficiently.

First consider the annulus $A_1$. Using the Imbalance Principle and the Twisting Number Lemma, we can increase $n_1$ by $1$, and repeat this process till either $n_1=0$ or $|2n_1-1|\leq |t|$, whichever comes first. If $n_1=0$ comes first, then we have $t(F)=-1$, which is a contradiction. This means that the procedure stops at an integer $n_1\leq-1$ with $|2n_1-1|\leq |t|$. If $|2n_1-1|<|t|$, then we can use the Imbalance Principle to increase the twisting number of $F$, which contradicts our choice of $F$. So $|2n_1-1|=|t|$, which implies that $t=2n_1-1\leq-3$.

Next consider the annulus $A_2$. Using the Imbalance Principle and the Twisting Number Lemma, we can increase $n_2$ by $1$, and repeat this process till either $n_2=-1$ or $|3n_2+1|\leq |t|$, whichever comes first. If $n_1=-1$ comes first, then we have $t(F)\geq -2$, which is a contradiction. This means that the procedure stops at an integer $n_2\leq-2$ with $|3n_2+1|\leq |t|$. If $|3n_2+1|<|t|$, then we can use the Imbalance Principle to increase the twisting number of $F$, which contradicts our choice of $F$. So $|3n_2+1|=|t|$, which implies that $t=3n_2+1\leq-5$.

Clearly, there is a positive integer $m$ satisfying $t=1-6m$, $n_1=1-3m$ and $n_2=-2m$. Now connect $T_1$ and $T_2$ by a vertical convex annulus $A$ with Legendrian boundary intersecting the dividing curves of $T_1,T_2$ efficiently. If $A$ has $\partial$-parallel dividing curves, then we can use the Legendrian Realization Principle to find a Legendrian regular fiber with twisting number greater than $t$, which contradicts our choice of $t$. So every dividing curve of $A$ connects one boundary component of $A$ to the other. Cut $\Sigma\times S^1$ along $A$ and smooth the edges. We get a torus $T$ isotopic to $T_3$ with dividing curves of slope $-\frac{m}{6m-1}$. If $m\geq n$, then
\[
-\frac{m}{6m-1} \geq -\frac{n}{6n-1} > s_3=-\frac{nn_3+1}{(6n-1)n_3+6},
\]
where $s_3=-\frac{nn_3+1}{(6n-1)n_3+6}$ is the slope of dividing curves of $T_3$. By \cite[Proposition 4.16]{H1}, there is a torus isotopic to $T_3$ with vertical dividing curves (isotopic to a regular fiber.) By the Legendrian Realization Principle, we can again find a Legendrian regular fiber with twisting number $0$, which is a contradiction. This shows that $m<n$.

The torus $T$ separates $-\Sigma(2,3,6n-1)$ into two sides. One side is a solid torus $V$ isotopic to $V_3$. The other side $(-\Sigma(2,3,6n-1))\setminus V$ is the union of $V_1$, $V_2$ and a neighborhood of the annulus $A$. The dividing curves of $A$ are unique up to an isotopy of $A$ fixing one boundary component since none of the dividing curves is $\partial$-parallel. Fix the dividing curves on $A$, since $V_1$ and $V_2$ are standard neighborhoods of Legendrian knots. it is easy to see that $\xi|_{(-\Sigma(2,3,6n-1))\setminus V}$ is uniquely determined up to isotopy relative to $T$. When measured in the coordinates of $V_3$, the slope of dividing curves of $T$ is $m-n$. So, by \cite[Theorem 2.3]{H1}, up to isotopy relative to $T$, there are $n-m$ tight contact structures on $V$ satisfying the given boundary condition. Note that, for each pair of possible dividing sets of $A$, there is an isotopy of $-\Sigma(2,3,6n-1)$ that maps one of them to the other. Thus, up to isotopy of $-\Sigma(2,3,6n-1)$, there are at most $n-m$ tight contact structures on $-\Sigma(2,3,6n-1)$ with twisting number $1-6m$. So the number of tight contact structures up to isotopy on $-\Sigma(2,3,6n-1)$ is at most
\[
\frac{n(n-1)}{2}=\sum_{m=1}^{n-1}(n-m).
\]
\end{proof}

It seems that the number of tight contact structures on $-\Sigma(2,3,6n-1)$ is exactly $\frac{n(n-1)}{2}$ since there are actually
$\frac{n(n-1)}{2}$ different Legendrian surgery constructions of
tight contact structures on $-\Sigma(2,3,6n-1)$. Before constructing
these surgeries, we need some preliminaries about tight contact structures on
the small Seifert fibered manifold
$M(-\frac{1}{2},\frac{1}{3},\frac{1}{6})$, which is also the torus
bundle over $S^1$ given by the monodromy map $\psi:T^2\rightarrow
T^2$ induced by
\[
\Psi = \left(
         \begin{array}{cc}
           1 & 1 \\
           -1 & 0 \\
         \end{array}
       \right):
\mathbb{R}^2\rightarrow\mathbb{R}^2.
\]

\begin{proposition}\cite[Theorem 0.1]{H2}\label{class-236}
There is a sequence of pairwise non-isotopic tight contact
structures $\{\xi_m\}_{m=1}^\infty$ on
$M(-\frac{1}{2},\frac{1}{3},\frac{1}{6})$. Any tight contact
structure on $M(-\frac{1}{2},\frac{1}{3},\frac{1}{6})$ is isotopic
to one of the $\xi_m$'s.
\end{proposition}

\begin{proposition}\cite[Propositions 15 and 16]{DG}\label{filling-236}
There is a simply connected symplectic manifold $(W,\omega)$ that
weakly fills $(M(-\frac{1}{2},\frac{1}{3},\frac{1}{6}),\xi_m)$ for
$\forall~ m\geq1$.
\end{proposition}
\begin{proof}
Such a symplectic manifold $(W,\omega)$ is constructed in \cite[Propositions 15 and 16]{DG}. We only need to show that $W$
is simply connected. Note that
\[
\left(
  \begin{array}{cc}
    1 & 1 \\
    -1 & 0 \\
  \end{array}
\right)
=
\left(
  \begin{array}{cc}
    1 & 0 \\
    -1 & 1 \\
  \end{array}
\right)
\left(
  \begin{array}{cc}
    1 & 1 \\
    0 & 1 \\
  \end{array}
\right).
\]
By the construction of $W$, there is a Lefschetz fibration
$W\rightarrow D^2$ which has exactly two singular points. The
vanishing circles of these two singular points induce a
$\mathbb{Z}$-basis for $\pi_1(T^2) \cong H_1(T^2) \cong \mathbb{Z}^2$. By \cite[Proposition 8.1.9]{GoSt}, there is an exact sequence
\[
\pi_1(T^2)\rightarrow\pi_1(W)\rightarrow\pi_1(D^2)(=0),
\]
where the first map is induced by the inclusion of $T^2$ into $W$ as
a regular fiber, and the second is induced by the projection. It follows that $\pi_1(W)=0$.
\end{proof}

The point $(0,0)^T\in\mathbb{R}^2$ induces the unique fixed point of
$\psi$, and gives a closed orbit $K_0$ in
$M(-\frac{1}{2},\frac{1}{3},\frac{1}{6})$, which is isotopic to the
$\frac{1}{6}$-singular fiber of the Seifert fibration. The torus
bundle structure gives $K_0$ a standard framing (c.f. \cite{Gh1}).
For any Legendrian knot $K$ in a tight contact manifold
$(M(-\frac{1}{2},\frac{1}{3},\frac{1}{6}),\xi)$ that is smoothly
isotopic to $K_0$, define its twisting number $t(K)$ to be the index
of its contact framing relative to this standard framing. Denote by
$t(\xi)$ the maximum of all such twisting numbers.

\begin{proposition}\cite[Lemma 3.5]{Gh1}\label{maxtwist-236}
$t(\xi_m)=-m$.
\end{proposition}

Performing a $(-n)$-surgery along the $\frac{1}{6}$-singular fiber
of $M(-\frac{1}{2},\frac{1}{3},\frac{1}{6})$ with respect to the
standard framing, we get $-\Sigma(2,3,6n-1)\cong
M(-\frac{1}{2},\frac{1}{3},\frac{n}{6n-1})$. For each
$m\in\{1,\cdots,n-1\}$, let $K_0^{(m)}$ be a Legendrian knot in
$(M(-\frac{1}{2},\frac{1}{3},\frac{1}{6}),\xi_m)$ smoothly
isotopic to $K_0$ with
$t(K_0^{(m)})=-m$. If we stabilize $K_0^{(m)}$ $n-m-1$ times, and
then perform a Legendrian surgery on resulted Legendrian knot, we
get a weakly fillable tight contact structure on
$-\Sigma(2,3,6n-1)$. (It is actually strongly fillable since
$-\Sigma(2,3,6n-1)$ is an integral homology sphere, c.f. \cite{E6}.)
There are $n-m$ ways to perform such an iterated stabilization of
$K_0^{(m)}$ depending on the number of positive stabilizations used
in the process. Denote by $\xi_{m,p}$ the tight contact structure on
$-\Sigma(2,3,6n-1)$ from the iterated stabilization of $K_0^{(m)}$
with $p$ positive stabilization and $n-m-p-1$ negative
stabilizations. This gives us a tight contact structure on
$-\Sigma(2,3,6n-1)$ for each pair $(m,p)$, where $1\leq m \leq n-1$,
and $0\leq p \leq n-m-1$. Altogether, we get $\frac{n(n-1)}{2}$
tight contact structures on $-\Sigma(2,3,6n-1)$. The hard part is to
show these tight contact structures are pairwise non-isotopic. Using
Theorem \ref{surgery-main}, we have following partial result.

\begin{proposition}\label{fix-m}
If $0\leq p_1,p_2 \leq n-m-1$ and $p_1\neq p_2$, then $\xi_{m,p_1}$
and $\xi_{m,p_2}$ are not isotopic.
\end{proposition}

\begin{proof}
This is a straightforward consequence of part (1) of Theorem
\ref{surgery-main} and the fact that the symplectic filling
$(W,\omega)$ of $\xi_m$ is simply connected.
\end{proof}

\begin{proof}[Proof of Theorem \ref{2n-3}]
Let $V$ be the cobordism from
$M(-\frac{1}{2},\frac{1}{3},\frac{1}{6})$ to $-\Sigma(2,3,6n-1)$
induced by the $(-n)$-surgery along the $\frac{1}{6}$-singular
fiber. Then, from Theorem \ref{OS-surgery}, we know that
$F^+_V(c^+(\xi_{m,p}))=c^+(\xi_m)$. Ghiggini \cite{Gh1} showed
that $c^+(\xi_1)\neq c^+(\xi_2)$. So $\xi_{1,p_1}$ and $\xi_{2,p_2}$
are non-isotopic for $0\leq p_1\leq n-2$ and $0\leq p_2\leq n-3$.
Combine this with Proposition \ref{fix-m}, we know that
$\xi_{1,0},\cdots,~\xi_{1,n-2},~\xi_{2,0},\cdots,~\xi_{2,n-3}$ are
$2n-3$ pairwise non-isotopic tight contact structures on
$-\Sigma(2,3,6n-1)$.
\end{proof}

The author hopes that, by a more careful computation of the
Ozsv\'ath-Szab\'o contact invariants, we can strengthen Theorem \ref{2n-3} and show that $\xi_{m_1,p_1}$ and $\xi_{m_2,p_2}$ are not
isotopic when $m_1\neq m_2$, which would complete the classification
of tight contact structures on $-\Sigma(2,3,6n-1)$.

\section{An example where our method does not apply}

The author was informed of Example \ref{counter} by Ghiggini,
which was proposed by Stipsicz.

\begin{example}\label{counter}
Consider the Stein fillable contact structure on $S^2\times S^1$.
Let $K$ be any Legendrian knot that is smoothly isotopic to an
$S^1$-fiber. Perform a Legendrian surgery on $K$, we get a Stein
fillable contact $3$-manifold, where the underlying smooth
$3$-manifold is $S^3$. To see this, note that $S^2\times S^1$ can
be constructed by performing a $0$-surgery on an unknot in $S^3$,
and an $S^1$-fiber comes from another unknot that links once with
the surgery unknot. So, topologically, the result of performing a
Legendrian surgery along $K$ is the same as performing a surgery
along a Hopf link in $S^3$, where one of its components has
coefficient $0$, and the other has an integer coefficient. This
clearly gives $S^3$. But there is only one tight contact structure
on $S^3$. This means the result of the Legendrian surgery here
does not depend on the choice of the Legendrian knot.
\end{example}

Ghiggini further remarked that, in the setting of Theorem
\ref{surgery-main}, if $[K]$ is a primitive element of $H_1(M)$,
then $H^2((M\times I) \cup_{\varphi_i} \mathcal{H}_2)=H^2(M)$, and
there is a unique $Spin^{\mathbb{C}}$-structure on $(M\times I)
\cup_{\varphi_i} \mathcal{H}_2$ that extends the
$Spin^{\mathbb{C}}$-structure on $M$ given by the contact structure.
So it is not possible to use $Spin^\mathbb{C}$-structures on the
cobordism to distinguish between contact structures resulted from the
Legendrian surgeries on stabilizations of $K$. Clearly, in the
weakly fillable case of Theorem \ref{surgery-main}, if $[K]$ is a
primitive element of $H_1(W)$, a similar remark applies. (These
examples correspond to the situation when $d=1$ in Theorem
\ref{surgery-main}. And Theorem \ref{surgery-main} does not give any
information about the result contact structures when $d=1,2$.)

\section*{Acknowledgments}
The author would like to thank Tomasz Mrowka for motivation and
helpful discussions, Paolo Ghiggini for pointing out a mistake in
an earlier version of this paper, and Andr\'as Stipsicz for
providing Example \ref{counter} to illustrate the mistake. The
author would also like to thank Peter Ozsv\'ath for helpful
discussions about \cite{OS1}.

\end{document}